\documentclass[twoside]{amsart}

\usepackage{latexsym, amsmath, amssymb, amsfonts, amscd, amsthm, mathrsfs}

\newcommand{\Ce}{{\mathbb C}}
\renewcommand{\Re}{{\mathbb R}}

\newcommand{\Ze}{{\mathbb Z}}
\newcommand{\Ne}{{\mathbb N}}
\newcommand{\Te}{{\mathbb T}}

\newcommand{\spn}{{\rm span}}

\newcommand{\pa}{\parallel}

\newcommand{\ie}{{\it i.e. }}

\theoremstyle{plain}
\newtheorem{theorem}{Theorem}[section]
\newtheorem{lemma}[theorem]{Lemma}
\newtheorem{proposition}[theorem]{Proposition}

\theoremstyle{definition}
\newtheorem{remark}[theorem]{Remark}

\newtheorem{notation}[theorem]{Notation}
\newtheorem{definition}[theorem]{Definition}

\newtheorem*{acknowledgments}{Acknowledgments}

\begin{document}

\title[Strong Morita Equivalence]
{Strong Morita Equivalence of Higher-dimensional Noncommutative Tori}
\author{Hanfeng Li}

\address{Department of Mathematics \\
University of Toronto \\
Toronto ON M5S 3G3, CANADA} \email{hli@fields.toronto.edu}
\date{December 5 , 2003}


\begin{abstract}
We show that matrices in the same orbit of the $SO(n, n|\Ze)$
action on the space of $n\times n$ skew-symmetric matrices give
strongly Morita equivalent noncommutative tori, both at the
$C^*$-algebra level and at the smooth algebra level. This proves a
conjecture of Rieffel and Schwarz.
\end{abstract}

\maketitle

\section{Introduction}
\label{intro:sec}

Let $n\ge 2$ and $\mathcal{T}_n$ be the space of $n\times n$ real
skew-symmetric matrices. For each $\theta \in \mathcal{T}_n$ the
corresponding $n$-dimensional noncommutative torus $A_{\theta}$ is
defined as the universal $C^*$-algebra generated by unitaries
$U_1, \cdots, U_n$ satisfying the relation
\begin{eqnarray*}
U_kU_j=e(\theta_{kj})U_jU_k,
\end{eqnarray*}
where $e(t)=e^{2\pi it}$. Noncommutative tori are one of the
canonical examples in noncommutative differential geometry
\cite{Rieffel90, Connes94}.

One may also consider the smooth version $A^{\infty}_{\theta}$ of
a noncommutative torus, which is the algebra of formal series
\begin{eqnarray*}
\sum c_{j_1, \cdots, j_n}U^{j_1}_1\cdots U^{j_n}_n,
\end{eqnarray*}
where the coefficient function  $\Ze^n\ni(j_1, \cdots, j_n)\mapsto
c_{j_1, \cdots, j_n}$ belongs to the Schwartz  space
$\mathcal{S}(\Ze^n)$  i.e. the space of $\Ce$-valued functions on $\Ze^n$
which vanish at infinity more
rapidly than any polynomial grows.
This is the space of smooth elements of
$A_{\theta}$ for the canonical action of $\mathbb{T}^n$ on
$A_{\theta}$.

The notion of (strong) Morita equivalence of $C^*$-algebras was
introduced by Rieffel \cite{Rieffel74, Rieffel82}. Strongly Morita
equivalent $C^*$-algebras share a lot of important properties such
as equivalent categories of modules, isomorphic $K$-groups, etc.,
and hence are usually thought to have the same geometry. In
\cite{Schwarz98} Schwarz also introduced the notion of complete
Morita equivalence of smooth noncommutative tori (see
Section~\ref{modules:sec} below), which is stronger than strong
Morita equivalence and has important application in M(atrix)
theory \cite{Schwarz98, KS02}.

A natural question is to classify noncommutative tori up to strong
Morita equivalence. Such results have important application to
physics  \cite{CDS97, Schwarz98}. For $n=2$ this was done by
Rieffel \cite{Rieffel81}. In this case there is a (densely
defined) action of the group $GL(2, \Ze)$ on $\mathcal{T}_2$, and
two matrices in $\mathcal{T}_2$ give strongly Morita equivalent
noncommutative tori if and only if they are in the same orbit of
this action. The higher dimensional case is much more complicated.
In \cite{Rieffel99a} Rieffel and Schwarz found a (densely defined)
action of $SO(n, n|\Ze)$ on $\mathcal{T}_n$ generalizing the above
$GL(2, \Ze)$-action. Here $O(n, n|\Re)$ is the group of linear
transformations of the space $\Re^{2n}$ preserving the quadratic
form $x_1x_{n+1}+x_2x_{n+2}+\cdots+x_nx_{2n}$, and $SO(n, n|\Ze)$
is the subgroup of $O(n, n|\Re)$ consisting of matrices with integer
entries and determinant $1$.

Following \cite{Rieffel99a} we write the elements of $O(n, n|\Re)$
in $2\times 2$ block form:
\begin{eqnarray*}
g=\begin{pmatrix} A & B \\ C &  D \end{pmatrix}.
\end{eqnarray*}
Here $A, B, C, D$ are $n\times n$ matrices satisfying
\begin{eqnarray} \label{O(n,n|R):eq}
A^tC+C^tA=0=B^tD+D^tB, & A^tD+C^tB=I.
\end{eqnarray}
The action of $SO(n, n|\Ze)$ is then defined as
\begin{eqnarray} \label{action:eq}
g\theta=(A\theta+B)(C\theta+D)^{-1},
\end{eqnarray}
whenever $C\theta+D$ is invertible. For each $g\in SO(n, n|\Ze)$
this action is defined on a dense open subset of $\mathcal{T}_n$.

Rieffel and Schwarz conjectured that if two matrices in
$\mathcal{T}_n$ are in the same orbit of this action then they
give strongly Morita equivalent noncommutative tori, both at the
$C^*$-algebra level and at the smooth algebra level. They proved
it for matrices restricted to a certain subset of $\mathcal{T}_n$
of second category. They also showed that the converse of their
conjecture at the $C^*$-algebra level fails for $n=3$
\cite[page 297]{Rieffel99a}, in contrast
to the case $n=2$, using the classification results of G. A.
Elliott and Q. Lin \cite{Lin96}.

The main goal of this paper is to prove their conjecture:

\begin{theorem} \label{orbit implies Morita equiv:theorem}
For any $\theta\in \mathcal{T}_n$ and $g\in SO(n, n|\Ze)$, if
$g\theta$ is defined then $A_{\theta}$ and $A_{g\theta}$ are
strongly Morita equivalent. Also $A^{\infty}_{\theta}$ and
$A^{\infty}_{g\theta}$ are completely Morita equivalent.
\end{theorem}

Schwarz has proved that if two matrices in $\mathcal{T}_n$ give
completely Morita equivalent smooth noncommutative tori then they
are in the same orbit of the $SO(n, n|\Ze)$-action
\cite[Section 5]{Schwarz98}. Thus we get

\begin{theorem} \label{complete Morita equiv:theorem}
Two matrices in $\mathcal{T}_n$ give completely Morita equivalent
smooth noncommutative tori if and only if they are in the same
orbit of the $SO(n, n|\Ze)$-action.
\end{theorem}

We have learned recently that using classification theory N. C.
Phillips has been able to show that two simple noncommutative tori
$A_{\theta}$ and $A_{\theta'}$ are strongly Morita equivalent if
and only if their ordered $K_0$-groups are isomorphic \cite[Remark
7.9]{Phillips03}. It would be interesting to see directly from the
matrices why the ordered $K_0$-groups of $A_{\theta}$ and
$A_{g\theta}$ are isomorphic.

This paper is organized as follows. Our proof of
Theorem~\ref{orbit implies Morita equiv:theorem} is constructive,
and we shall use the Heisenberg equivalence modules constructed by
Rieffel in \cite{Rieffel88}. So we recall briefly Rieffel's
construction first in Section~\ref{modules:sec}. In order to apply
Rieffel's construction we need to reduce an arbitrary matrix in
$\mathcal{T}_n$ to one satisfying certain nice properties. This is
done in Section~\ref{matrices:sec}. We prove Theorem~\ref{orbit
implies Morita equiv:theorem} in Section~\ref{proof:sec}.

\begin{acknowledgments}
  I would like to thank Marc Rieffel for many helpful discussions and suggestions.
  I also thank Albert Schwarz for valuable discussions about complete Morita equivalence.
  I also thank Rolf Svegstrup for pointing out some misprints and the referee for helpful
comments making the paper more readable.
\end{acknowledgments}

\section{Heisenberg Equivalence Modules}
\label{modules:sec}

In this section we recall Schwarz's definition of complete Morita
equivalence and Rieffel's construction of Heisenberg equivalence
modules for noncommutative tori.

Let $L=\Re^n$. We shall think of $\Ze^n$ as the standard lattice
in $L^*$, and $\theta$ as in $\wedge^2L$. One may also describe
$A_{\theta}$ as the universal $C^*$-algebra generated by unitaries
$\{U_x\}_{x\in \Ze^n}$ satisfying the relation
\begin{eqnarray} \label{tori:eq}
U_xU_y=\sigma_{\theta}(x, y)U_{x+y},
\end{eqnarray}
where we write $x, y$ as column vectors, and
$\sigma_{\theta}(x, y)= e((x\cdot \theta y)/2)$. Under this
description the smooth algebra $A^{\infty}_{\theta}$ becomes
$\mathcal{S}(\Ze^n, \sigma_{\theta})$, the Schwartz space
$\mathcal{S}(\Ze^n)$ equipped with the convolution induced by
(\ref{tori:eq}). There is a canonical action of the Lie algebra
$L$ as derivations on $A^{\infty}_{\theta}$, which is induced by
the canonical action of $\mathbb{T}^n$ on $A_{\theta}$ and is
given explicitly by
\begin{eqnarray*}
\delta_X(U_x)=2\pi i \left<X, x\right>U_x
\end{eqnarray*}
for all $X\in L$ and $x\in \Ze^n$, where $\left<\cdot, \cdot\right>$ denotes
the natural pairing between $L$ and $L^*$.

Given a right $A^{\infty}_{\theta}$-module $E$, a {\it connection
} on $E$ is a linear map $\nabla:L\rightarrow Hom_{\Ce}(E)$
satisfying the Leibniz rule:
\begin{eqnarray*}
\nabla_X(fU_x)=(\nabla_Xf)U_x+f\cdot \delta_X(U_x)
\end{eqnarray*}
for all $X\in L, \, f\in E$ and $x\in \Ze^n$.
For each $X\in L$ the connection $\nabla$ induces a derivation
$\hat{\delta}_X$ on $End_{A^{\infty}_{\theta}}(E)$ by
\begin{eqnarray*}
(\hat{\delta}_Xa)(f)=\nabla_X(a f)-a\cdot \nabla_X f
\end{eqnarray*}
for all $a\in End_{A^{\infty}_{\theta}}(E)$ and $f\in E$.
If $\nabla$ has {\it constant curvature}, i.e. there is
skew-symmetric bilinear map $\Omega:L\times L\rightarrow \Ce$ such
that $[\nabla_X, \nabla_Y]=\Omega(X, Y)\cdot 1$ for all $X,\, Y\in L$,
then $X\mapsto
\hat{\delta}_X$ is a Lie algebra homomorphism from $L$ to the
derivation space $Der(End_{A^{\infty}_{\theta}}(E))$ of
$End_{A^{\infty}_{\theta}}(E)$. When $E$ is equipped with an
$A^{\infty}_{\theta}$-valued inner product, we shall consider only
{\it Hermitian } connections, i.e. $\delta_X\left<f,
g\right>=\left<\nabla_X f,\, g\right>+\left<f,\, \nabla_X g\right>$ for $X\in L$ and
$f, \, g\in A^{\infty}_{\theta}$.

 We refer to \cite{Rieffel82} for the
definition and standard facts about strong Morita equivalence of
$C^*$-algebras. Let $E$ be a strong Morita equivalence
$A^{\infty}_{\theta'}$-$A^{\infty}_{\theta}$-bimodule. For clarity
we let $L_{\theta}$ and $L_{\theta'}$ denote the space $L$ for
$\theta$ and $\theta'$ respectively. We say that $E$ is a {\it
complete Morita equivalence
$A^{\infty}_{\theta'}$-$A^{\infty}_{\theta}$-bimodule}
\cite[page 729]{Schwarz98} if there is a constant-curvature connection
$\nabla$ on $E_{A^{\infty}_{\theta}}$ and a linear isomorphism
$\phi:L_{\theta}\rightarrow L_{\theta'}$ such that the induced Lie
algebra homomorphism $L_{\theta}\rightarrow
Der(End_{A^{\infty}_{\theta}}(E))=Der(A^{\infty}_{\theta'})$
coincides with the composition homomorphism
$L_{\theta}\overset{\phi}\rightarrow L_{\theta'}\rightarrow
Der(A^{\infty}_{\theta'})$. Intuitively, this means that the
equivalence bimodule $E$ is "smooth", i.e. it transfers the
tangent spaces ($L_{\theta}$ and $L_{\theta'}$) of the
noncommutative differentiable manifolds $A^{\infty}_{\theta}$ and
$A^{\infty}_{\theta'}$ back and forth.

Next we recall Rieffel's construction of Heisenberg equivalence
modules in \cite[Sections 2-4]{Rieffel88}. Let $M$ be a locally
compact abelian group, let $\hat{M}$ be its dual group, and let
$G=M\times \hat{M}$. There is a canonical Heisenberg cocycle on
$G$ defined by
\begin{eqnarray*}
\beta((m,s),(l,t))=\left<m, t\right>,
\end{eqnarray*}
where $\left<\cdot, \cdot\right>$ denotes the natural pairing
between $M$ and $\hat{M}$. There is also a skew bicharacter,
$\rho$, on G defined by
\begin{eqnarray} \label{rho:eq}
\rho(x,y)=\beta(x,y)\bar{\beta}(y,x).
\end{eqnarray}

We'll concentrate on the case $M=\Re^p\times \Ze^q\times W$, where
$p, q\in \Ze_{\ge 0}$ with $2p+q=n$ and $W$ is a finite abelian
group. Say $W=\Ze_{n_1}\times \dots \times \Ze_{n_k}$ for some
$n_1, \cdots, n_k\in \Ne$. We shall write $G$ as $\Re^{p}\times
\Re^{*p}\times \Ze^q\times \Te^q\times (\Ze_{n_1}\times \dots
\times \Ze_{n_k})\times (\Ze_{n_1}\times \dots \times \Ze_{n_k})$.
Let
\begin{eqnarray} \label{J:eq}
P_1=diag(\frac{1}{n_1}, \cdots, \frac{1}{n_k}), & & J_0=\begin{pmatrix} 0 & I_p \\ -I_p & 0
\end{pmatrix}, \quad J_1=\begin{pmatrix} J_0 & 0 & 0\\ 0 &  0 & I_q \\
0 & -I_q & 0\end{pmatrix},\\
\nonumber
J_2=\begin{pmatrix} 0 & P_1\\ -P_1 &0
\end{pmatrix}, & &\,  J=\begin{pmatrix} J_1 & 0 \\ 0 & J_2\end{pmatrix}.
\end{eqnarray}
Then $J$ is a square matrix of size $n+q+2k$, and we shall think
of it as a $2$-form on $H^*:=\Re^{p}\times \Re^{*p}\times
\Re^q\times \Re^{*q}\times \Re^k\times \Re^{*k}$. Let $J'$ be the
matrix obtained by replacing negative entries of $J$ by $0$. Then
$J=J'-(J')^t$. For any $x, y\in G$
we have
\begin{eqnarray*}
\beta(x,y)=e(x\cdot J'y) \quad \mbox{ and } \quad \rho(x,y)=e(x\cdot J y),
\end{eqnarray*}
where we
use the natural covering map $\Re^{p}\times
\Re^{*p}\times \Ze^q\times \Re^{*q}\times \Ze^k\times \Ze^k\rightarrow G$
to write $x$ and $y$ as column vectors
in $\Re^{n+q+2k}$ (notice that though $J'y$ depends on the choice
of the representative of $y$ in $\Re^{p}\times
\Re^{*p}\times \Ze^q\times \Re^{*q}\times \Ze^k\times \Ze^k$,
the
values $e(x \cdot J'y)$ and $e(x\cdot J y)$ do not depend on such choice).

\begin{definition}\cite[Definition 4.1]{Rieffel88} \label{embedding map:def}
By an
\emph{embedding map} for $\theta \in \mathcal{T}_n$ we mean a
linear map $T$ from $L^*$ to $H^*$ such that:

(1) $T(\Ze^n)\subseteq \Re^{p}\times \Re^{*p}\times \Ze^q\times
\Re^{*q}\times \Ze^k\times \Ze^k$. Then we can think of $T(\Ze^n)$
as in $G$ via composing $T|_{\Ze^n}$ with the natural covering map
$\Re^{p}\times \Re^{*p}\times \Ze^q\times \Re^{*q}\times
\Ze^k\times \Ze^k \rightarrow G$.

(2) $T(\Ze^n)$ is a lattice in $G$.

(3) The form $J$ on $H^*$ is pulled back by $T$ to the form
$\theta$ on $L^*$, i.e. $T^tJT=\theta$.
\end{definition}

The condition (2) above is equivalent to

(2') The map $\tilde{T}:=\gamma \circ T:L^*\rightarrow
\Re^{p}\times \Re^{*p}\times \Re^q$ is invertible, where $\gamma$
is the projection of $H^*$ onto $\Re^{p}\times \Re^{*p}\times
\Re^q$.

The bimodule Rieffel constructed is the Schwartz space
$\mathcal{S}(M)$, i.e. the space of smooth functions on $M$ which,
together with all their derivatives, vanish at infinity more
rapidly than any polynomial grows.

\begin{proposition}\cite[Theorem 2.15, Corollary 3.8]{Rieffel88} \label{module:prop}
Let $\theta, \theta'\in \mathcal{T}_n$, and let $T, S$ be
embedding maps of $L^*$ into $H^*$ for $\theta$ and $-\theta'$
respectively such that $S(\Ze^n)=(T(\Ze^n))^{\perp}$, where
$(T(\Ze^n))^{\perp}=\{z\in G: \rho(z, y)=1 \mbox{ for all } y\in
T(\Ze^n)\}$. Let $T'$ and $T''$ be the composition maps
$\Ze^n\overset{T}\rightarrow G\rightarrow M$ and
$\Ze^n\overset{T}\rightarrow G\rightarrow \hat{M}$ respectively.
Define $S'$ and $S''$ similarly. Fix a Haar measure on $M$. Then
$\mathcal{S}(M)$ is a strong Morita equivalence
$A^{\infty}_{\theta'}$-$A^{\infty}_{\theta}$-bimodule with the
module structure and inner products defined by:
\begin{eqnarray*}
(fU_x)(m)&=&e(-T(x)\cdot J'T(x)/2)\left<m,
T''(x)\right>f(m-T'(x)),\\
\left<f, g\right>_{\mathcal{S}(\Ze^n,
\sigma_{\theta})}(x)&=&e(-T(x)\cdot J'T(x)/2)\int_{G}\left<m,
-T''(x)\right>g(m+T'(x))\bar{f}(m)\, dm,\\
(V_xf)(m)&=&e(-S(x)\cdot J'S(x)/2)\left<m,
-S''(x)\right>f(m+S'(x)), \\
{}_{\mathcal{S}(\Ze^n, \sigma_{\theta'})}\left<f, g\right>(x)
&=&K\cdot e(S(x)\cdot J'S(x)/2)\int_{G}\left<m,
S''(x)\right>f(m)\bar{g}(m+S'(x))\, dm,
\end{eqnarray*}
where $K$ is a positive constant and for clarity $V_x$ denotes the
unitary in $\mathcal{S}(\Ze^n, \sigma_{\theta'})$. Moreover, there
is a linear map $Q:\Re^{*p}\times \Re^p\times \Re^{*q}\rightarrow
Hom_{\Ce}(\mathcal{S}(M))$ such that
$\nabla_X=Q_{(\tilde{T}^{-1})^*(X)}$ and
$\nabla'_X=Q_{(\tilde{S}^{-1})^*(-X)}$ are connections with
respect to $\mathcal{S}(M)_{A^{\infty}_{\theta}}$ and
${}_{A^{\infty}_{\theta'}}\mathcal{S}(M)$ respectively. The
connection $\nabla$ has constant curvature
\begin{eqnarray*}
\Omega=2\pi i\tilde{T}^{-1}(\sum^p_{j=1} \bar{e}_j\wedge e_j),
\end{eqnarray*}
where $e_1, \cdots, e_p$ are the standard basis of $\Re^p$
and $\bar{e}_1, \cdots, \bar{e}_p$ are the dual basis of
$\Re^{*p}$. Thus $\mathcal{S}(M)$ is a complete Morita equivalence
$A^{\infty}_{\theta'}$-$A^{\infty}_{\theta}$-bimodule. When
completed with the norm $\pa f\pa :=\pa \left< f,
f\right>_{A^{\infty}_{\theta}}\pa^{\frac{1}{2}}=\pa
{}_{A^{\infty}_{\theta'}}\left< f, f\right>\pa^{\frac{1}{2}}$,
$\mathcal{S}(M)$ becomes a strong Morita equivalence
$A_{\theta'}$-$A_{\theta}$-bimodule.
\end{proposition}

\begin{remark} \label{module:remark}
(1) The definition of embedding maps in Definition~\ref{embedding
map:def} differs from that in \cite[Definition 4.1]{Rieffel88} by
a sign of $\theta$. This is because Rieffel's $A_{\theta}$ is our
$A_{-\theta}$ (see the discussion at the end of page 285 of
\cite{Rieffel88}).

(2) In \cite[Section 4]{Rieffel88} the definition of embedding
maps and the part of Proposition~\ref{module:prop} above
concerning connections and curvatures are only given for the case
$W=0$ \cite[Definition 4.1]{Rieffel88} \cite[pages
290-291]{Rieffel88}. The general case was discussed there in terms
of tensor products with finite dimensional representations
\cite[Section 5]{Rieffel88}. For our purpose it's better to deal
with $\Re^p\times \Ze^q\times W$ directly. The proofs in
\cite[pages 290-291]{Rieffel88} for the case $W=0$ are easily
checked to hold for the general case.
\end{remark}
\section{Decomposition of Matrices}
\label{matrices:sec}

In Proposition~\ref{embedding:prop} we shall use the construction
in \cite{AS01} to find the appropriate finite abelian group $W$.
To this goal we need the matrix $g\in SO(n, n|\Ze)$ to be of the
special form in Lemma~\ref{exist H:lemma} below. We shall prove in
Lemma~\ref{normalize:lemma} that every $g$ can be reduced to such
a special one.

\begin{lemma} \label{DC^t ss:lemma}
Let $g=\begin{pmatrix} A & B \\ C &  D \end{pmatrix}\in
O(n,n|\Re)$. Then $DC^t$ is skew-symmetric.
\end{lemma}
\begin{proof}
Since $g\in O(n,n|\Re)$ we have that
\begin{eqnarray*}
\begin{pmatrix} A & B \\ C &  D \end{pmatrix}^t
\begin{pmatrix} 0 & I \\ I &  0 \end{pmatrix}
\begin{pmatrix} A & B \\ C &  D \end{pmatrix}
=\begin{pmatrix} 0 & I \\ I &  0 \end{pmatrix}.
\end{eqnarray*}
Hence
\begin{eqnarray*}
g^{-1}=
\begin{pmatrix} A & B \\ C &  D \end{pmatrix}^{-1}
=\begin{pmatrix} 0 & I \\ I &  0 \end{pmatrix}^{-1}
\begin{pmatrix} A & B \\ C &  D \end{pmatrix}^t
\begin{pmatrix} 0 & I \\ I &  0 \end{pmatrix}
=\begin{pmatrix} D^t & B^t \\ C^t &  A^t \end{pmatrix} .
\end{eqnarray*}
Since $O(n,n|\Re)$ is a group we have that $g^{-1}\in O(n,n|\Re)$.
By (\ref{O(n,n|R):eq}) the matrix $(D^t)^tC^t=DC^t$ is skew-symmetric.
\end{proof}

Using Lemma~\ref{DC^t ss:lemma} simple calculations yield:

\begin{lemma} \label{(C theta +D)^-1 C ss:lemma}
Let $g=\begin{pmatrix} A & B \\ C &  D \end{pmatrix}\in
O(n,n|\Re)$. Let $\theta \in \mathcal{T}_n$ with $C\theta +D$
invertible. Then $(C\theta +D)^{-1}C$ is skew-symmetric.
\end{lemma}

\begin{lemma} \label{exist H:lemma}
Let $g=\begin{pmatrix} A & B \\ C &  D
\end{pmatrix}\in SO(n,n|\Ze)$, and let $p\in \Ze_{\ge 0}$. Then
the following are equivalent:

(i) there is some $\theta \in \mathcal{T}_n$ such that $(C\theta
+D)^{-1}C$ is of the form $\begin{pmatrix} F_{11} & 0 \\ 0 & 0
\end{pmatrix}$ for some $F_{11}\in GL(2p | \Re)$;

(ii) there exists a $Z\in \mathcal{T}_{2p}$ such that
\begin{eqnarray*}
C=\begin{pmatrix} C_{11} & 0 \\ C_{21} &  0 \end{pmatrix} \mbox{
and } \, D=\begin{pmatrix} -C_{11}Z & D_{12} \\ -C_{21}Z &
D_{22}\end{pmatrix},
\end{eqnarray*}
where $C_{11}\in M_{2p}(\Ze)$.

In this event, the matrix $\begin{pmatrix} C_{11} & D_{12} \\
C_{21} & D_{22} \end{pmatrix}$ is invertible. The matrix $Z$ is
unique, and its entries are all rational numbers. Also for any
$\theta'\in \mathcal{T}_n$ in the
block form $\begin{pmatrix} \theta'_{11} & \theta'_{12} \\
\theta'_{21} & \theta'_{22}\end{pmatrix}$, where $\theta'_{11}$
has size $2p\times 2p$, the matrix $C\theta'+D$ is invertible if
and only if $\theta'_{11}-Z$ is invertible. In this case
\begin{eqnarray} \label{C theta+D inverse:eq}
(C\theta'+D)^{-1}C=\begin{pmatrix} (\theta'_{11}-Z)^{-1} & 0 \\ 0
& 0 \end{pmatrix}.
\end{eqnarray}
\end{lemma}
\begin{proof}
(i)$\Rightarrow$(ii). From the assumption we have
$C\begin{pmatrix} I_{2p} & 0 \\ 0 &  0\end{pmatrix}=C$. Thus $C$
has the desired form in (ii). Notice that
\begin{eqnarray*}
\begin{pmatrix} C_{11} & 0 \\ C_{21} &  0 \end{pmatrix} = C
= (C\theta+D)\begin{pmatrix} F_{11} & 0 \\ 0 &  0 \end{pmatrix}
= \begin{pmatrix} (C_{11}\theta_{11}+D_{11})F_{11} & 0 \\
(C_{21}\theta_{11}+D_{21})F_{11}& 0
\end{pmatrix},
\end{eqnarray*}
where we are writing both $\theta$ and $D$ in block forms. Thus
$C_{j1}=(C_{j1}\theta_{11}+D_{j1})F_{11}$ for $j=1,2$. Let
$Z=\theta_{11}-(F_{11})^{-1}$. Then $D_{j1}=-C_{j1}Z$. By
Lemma~\ref{(C theta +D)^-1 C ss:lemma} the matrix $F_{11}$ is
skew-symmetric. Then so is $Z$.

(ii)$\Rightarrow$(i). For any $\theta' \in \mathcal{T}_n$ we have
\begin{eqnarray*}
C\theta'+D = \begin{pmatrix} C_{11} & D_{12}\\ C_{21} &  D_{22}
\end{pmatrix}
\begin{pmatrix} \theta'_{11}-Z & \theta'_{12} \\ 0 &  I
\end{pmatrix}.
\end{eqnarray*}
Take $\theta\in \mathcal{T}_n$ such that $C\theta+D$ is
invertible.
Then $\begin{pmatrix} C_{11} & D_{12} \\
C_{21} & D_{22} \end{pmatrix}$ is invertible. Therefore
$C\theta'+D$ is invertible if and only if $\theta'_{11}-Z$ is
invertible. In this case simple computations yield (\ref{C theta+D
inverse:eq}). In particular $(C\theta+D)^{-1}C$ has the form
described in (i). By varying $\theta$ slightly we may assume that
$\theta$ is rational, \ie the entries of $\theta$ are all rational numbers.
Then so are $F_{11}$ and
$Z=\theta_{11}-(F_{11})^{-1}$.
\end{proof}

\begin{notation} \label{rho:notation}
For any $R\in GL(n|\Ze)$ we denote by $\rho(R)$ the matrix
$\begin{pmatrix} R & 0 \\ 0 &  (R^{-1})^t \end{pmatrix}$ in $SO(n,
n|\Ze)$. For any $N\in \mathcal{T}_n\cap M_n(\Ze)$ we denote by
$\mu(N)$ the matrix $\begin{pmatrix} I & N \\ 0 &  I
\end{pmatrix}$ in $SO(n, n|\Ze)$.
\end{notation}

Notice that the noncommutative tori corresponding to the matrices
$\rho(R)\theta=R\theta R^t$ and $\mu(N)\theta=\theta+N$ are both
isomorphic to $A_{\theta}$.

\begin{lemma} \label{normalize:lemma}
Let $g=\begin{pmatrix} A & B \\ C &  D \end{pmatrix}$ in
$SO(n,n|\Ze)$. Then there exists an $R\in GL(n| \Ze)$ such that
$g\cdot \rho(R)$ satisfies the condition (1) in Lemma~\ref{exist
H:lemma} for some $p\in \Ze_{\ge 0}$.
\end{lemma}
\begin{proof}
Let $V=\{X\in \Re^n|CX=0\}$, and let $K=V\cap \Ze^n$. Since the
entries of $C$ are all integers, $K$ spans $V$. By the elementary
divisors theorem \cite[page 153, Theorem III.7.8]{Lang02} we can
find a basis $\beta_1, \cdots, \beta_n$ of $\Ze^n$, some integer
$1\le k\le n$ and positive integers $c_k, \cdots, c_n$ such that
$K$ is generated by $c_k\beta_k, \cdots, c_n\beta_n$. Then
$V=\spn(\beta_k, \cdots, \beta_n)$. Let $e_1, \cdots, e_n$ be the
standard basis of $\Ze^n$. Then $(\beta_1, \cdots,
\beta_n)=(e_1,\cdots, e_n)R$ for some $R\in GL(n| \Ze)$. Let
\begin{eqnarray*}
\begin{pmatrix} A' & B' \\ C' &  D' \end{pmatrix}=\begin{pmatrix} A & B \\ C & D
\end{pmatrix}\rho(R)\in SO(n,n|\Ze).
\end{eqnarray*}
Choose $\theta \in \mathcal{T}_n$ such that $C\theta+D$ is
invertible. Let $\theta'=R^{-1}\theta (R^{-1})^t\in
\mathcal{T}_n$. Now we need
\begin{lemma} \label{of special form:lemma}
$(C'\theta' +D')^{-1}C'$ is of the form $\begin{pmatrix} F_{11} & 0 \\
0 & 0
\end{pmatrix}$ for some $F_{11}\in GL(k-1 | \Re)$.
\end{lemma}
\begin{proof} In view of Lemma~\ref{(C theta +D)^-1 C ss:lemma} this is clearly
equivalent to saying that the vectors $X=(x_1, \cdots, x_n)^t$ in
$\Re^n$ satisfying $(C'\theta' +D')^{-1}C'X=0$ are exactly those
with $x_1=\cdots=x_{k-1}=0$. Notice that $
(C'\theta'+D')^{-1}C'=R^t(C\theta+D)^{-1}CR$. Hence $(C'\theta'
+D')^{-1}C'X=0$ if and only if $CRX=0$, if and only if $RX\in V$,
if and only if $(\beta_1, \cdots, \beta_n)X\in V$, if and only if
$x_1=\cdots=x_{k-1}=0$.
\end{proof}
Back to the proof of Lemma~\ref{normalize:lemma}.
By Lemma~\ref{(C theta +D)^-1 C ss:lemma} the matrix $F_{11}$ is
skew-symmetric. Since $F_{11}\in GL(k-1 | \Ze)$ we see that $k-1$
is even. This completes the proof of Lemma~\ref{normalize:lemma}.
\end{proof}

\section{Strong Morita Equivalence}
\label{proof:sec}

In this section we prove Theorem~\ref{orbit implies Morita
equiv:theorem}. We shall employ the notation in
Section~\ref{modules:sec} and Lemma~\ref{exist H:lemma}. In view
of Proposition~\ref{module:prop} the key is to find embedding
maps. This is established in the following

\begin{proposition} \label{embedding:prop}
Let $g=\begin{pmatrix} A & B \\ C &  D \end{pmatrix}$ in
$SO(n,n|\Ze)$ satisfying the conditions (1) and (2) in
Lemma~\ref{exist H:lemma} for some $p\in \Ze_{\ge 0}$. Then there
exist an $N\in \mathcal{T}_n\cap M_n(\Ze)$, an $R\in GL(n | \Ze)$,
a $g'\in SO(n,n|\Ze)$ and a finite abelian group $W$ such that
$g=\nu(N)\rho(R)g'$ and for any $\theta \in \mathcal{T}_n$ with
$C\theta+D$ invertible there are embedding maps $T,
S:L^*\rightarrow H^*$ for $\theta$ and $-g'\theta$ respectively
satisfying $S(\Ze^n)=(T(\Ze^n))^{\perp}$ (see
Definition~\ref{embedding map:def}(1) and Proposition~\ref{module:prop}
for the meaning of $(T(\Ze^n))^{\perp}$).
\end{proposition}
\begin{proof}
Let $Z$ be as in Lemma~\ref{exist H:lemma} for $g$. Then $Z$ is
rational, and hence there is some $m\in \Ne$ such that $mZ$ is
integral. Thinking of $mZ$ as a bilinear alternating form on
$\Ze^n$,
by \cite[page 598, Exercise XV.17]{Lang02} we can find an
$R\in GL(2p | \Ze)$, some integer $1\le k\le p$ and integers $h_1,
\cdots, h_k$ such that
\begin{eqnarray*}
mZ=R^t\begin{pmatrix} 0 & P  & 0 \\ -P &  0 & 0\\ 0 & 0 &
0\end{pmatrix}R,
\end{eqnarray*}
where $P=diag(h_1, \cdots, h_k)$. Let $m_j/n_j=h_j/m$ with
$(m_j,n_j)=1$ and $n_j>0$ for each $1\le j\le k$. Set
$W=\Ze_{n_1}\times \cdots \times \Ze_{n_k}$.

Let $\theta \in \mathcal{T}_n$ with
$C\theta+D$ invertible.
We are ready to construct an embedding map for $\theta$ now. Our
method is similar to that in the proof of the proposition in
\cite{Rieffel99a}. But our situation is more complicated since we
have to deal with the torsion part $W$.
Write $\theta$ in block form as in Lemma~\ref{exist H:lemma}.
By Lemmas~\ref{exist H:lemma} and \ref{(C theta +D)^-1 C ss:lemma}
the matrix $\theta_{11}-Z$ is invertible and skew-symmetric.
So we
can find a $T_{11}\in GL(2p | \Re)$ such that
$T^t_{11}J_0T_{11}=\theta_{11}-Z$, where $J_0$ is defined in (\ref{J:eq}).
Let $T_{31}=\theta_{21}$, and
let $T_{32}$ be any $q\times q$ matrix such that
$T_{32}-T^t_{32}=\theta_{22}$, where $q=n-2p$.
Let $P_2=diag(m_1, \cdots, m_k)$,
and let
\begin{eqnarray*}
T_1=\begin{pmatrix} T_{11} & 0  \\ 0 & I_q\\ T_{31} &
T_{32}\end{pmatrix},\, T_2=\begin{pmatrix} P_2 & 0 & 0\\ 0 &  I_k
& 0\end{pmatrix}\begin{pmatrix} R  & 0 \\ 0 & I_q \end{pmatrix},
\, T=\begin{pmatrix} T_1  \\ T_2\end{pmatrix}.
\end{eqnarray*}
Then $T_1, T_2$ and $T$ have sizes $(n+q)\times n, \, 2k\times n$
and $(n+q+2k)\times n$ respectively. Simple calculations yield
$T^tJT=\theta$, where $J$ is defined in (\ref{J:eq}).
Notice that as a linear map from $L^*=\Re^{*n}$ to $H^*=\Re^{p}\times \Re^{*p}\times
\Re^q\times \Re^{*q}\times \Re^k\times \Re^{*k}$,
$T$ carries the lattice $\Ze^n=\Ze^{2p}\times \Ze^q$ into $\Re^{p}\times
\Re^{*p}\times \Ze^q\times \Re^{*q}\times \Ze^k\times \Ze^k$. Also
observe that $\tilde{T}$ (see Definition~\ref{embedding map:def}(2'))
is given by the invertible matrix
 $\begin{pmatrix} T_{11} & 0  \\
0 & I_q \end{pmatrix}$. Thus the conditions in Definition~\ref{embedding map:def}
are satisfied and hence $T$ is an embedding map for $\theta$.

Let $\mathcal{D}=T(\Ze^n)$.
By Definition~\ref{embedding map:def}(1) we may think of $\mathcal{D}$
as in $G=\Re^{p}\times
\Re^{*p}\times \Ze^q\times \Te^q\times (\Ze_{n_1}\times \dots
\times \Ze_{n_k})\times (\Ze_{n_1}\times \dots \times \Ze_{n_k})$.
We need to find some embedding map of
$\Ze^n$ into $G$
with image being exactly $\mathcal{D}^{\perp}=
\{z\in G: \rho(z, y)=1 \mbox{ for all } y\in \mathcal{D}\}$,
where $\rho$ is defined in (\ref{rho:eq}).

For any $x\in G$, it is in $\mathcal{D}^{\perp}$ exactly if
$x\cdot JTz\in \Ze$ for all $z\in \Ze^n$, exactly if $T^tJx\in
\Ze^n$. Let $T_3=\begin{pmatrix} 0 \\ -I_q
\end{pmatrix}$
be an $(n+q)\times q$ matrix. Let $T_4=diag(n_1, \cdots, n_k, n_1,
\dots, n_k)$. Set
\begin{eqnarray*}
\bar{T}=\begin{pmatrix} T_1 & T_3 & 0\\ T_2 & 0 &T_4\end{pmatrix},
\end{eqnarray*}
a square matrix of size $n+q+2k$. It is easy to check that
$T^tJx\in \Ze^n$ exactly if $\bar{T}^tJx\in \Ze^{n+q+2k}$. Also it
is easy to see that $\bar{T}$ is invertible. Thus
\begin{eqnarray*}
\mathcal{D}^{\perp}=(\bar{T}^tJ)^{-1}(\Ze^{n+q+2k}).
\end{eqnarray*}
Recall the matrices $J_0$ and $J_1$ defined in (\ref{J:eq}).
Straight-forward  calculations show that
\begin{eqnarray*}
(\bar{T}^tJ)^{-1} = \begin{pmatrix}  -J_1\begin{pmatrix}  T^t_1\\
T^t_3
\end{pmatrix}^{-1} & J_1\begin{pmatrix}  T^t_1 \\ T^t_3 \end{pmatrix}^{-1}
\begin{pmatrix}T^t_2\\  0 \end{pmatrix}T^{-1}_4\\0 & \begin{pmatrix} 0 & -I_k \\ I_k & 0 \end{pmatrix}
\end{pmatrix},
\end{eqnarray*}
and
\begin{eqnarray*}
J_1\begin{pmatrix}  T^t_1\\ T^t_3
\end{pmatrix}^{-1}
=\begin{pmatrix} J_0(T^t_{11})^{-1} & 0 &
-J_0(T^t_{11})^{-1}T^t_{31} \\0 & 0 & I_q\\
 0 & -I_q & T^t_{32}\end{pmatrix}.
\end{eqnarray*}
Thus
\begin{eqnarray*}
(\bar{T}^tJ)^{-1}(0^{2p}\times \Ze^q\times 0^q\times
0^{2k})=0^{2p}\times 0^q\times \Ze^q\times 0^{2k},
\end{eqnarray*}
which is $0$ in $G$. So
$(\bar{T}^tJ)^{-1}(\Ze^{n+q+2k})=(\bar{T}^tJ)^{-1}(\Ze^{2p}\times
0^q\times \Ze^{q+2k})$. Let $\triangle$ be the set of all vectors
$y=\begin{pmatrix} y_1\\ 0 \\ y_2\\y_3 \end{pmatrix}$ in
$\Ze^{2p}\times 0^q\times \Ze^q\times \Ze^{2k}$ satisfying $y_3\in
n_1\Ze\times \cdots \times n_k\Ze\times n_1\Ze\times \cdots \times
n_k\Ze$ and $\begin{pmatrix} y_1\\ 0 \\ y_2
\end{pmatrix}=\begin{pmatrix}T^t_2\\  0 \end{pmatrix}T^{-1}_4y_3$.
Then $(\bar{T}^tJ)^{-1}(\triangle)=0$ in $G$. For each $1\le j\le
k$ since $(m_j, n_j)=1$ we can find $c_j, d_j\in \Ze$ such that
$c_jm_j+d_jn_j=1$. Now we need
\begin{lemma} \label{lifting:lemma}
Let $\varphi_1$ be the embedding $\Ze^n\hookrightarrow
\Ze^{2p}\times 0^q\times \Ze^{q+2k}$ sending $(x_1, \cdots, x_n)^t$
to
\begin{eqnarray*}
(-d_1x_1, \cdots, -d_kx_k,0, \cdots, 0, x_{2k+1}, \cdots,
x_{2p})^t\times 0^q\times \\
 (x_{2p+1}, \cdots, x_n, c_1x_1, \cdots,
c_kx_k, x_{k+1}, \dots, x_{2k})^t.
\end{eqnarray*}
Let $\varphi$ be the composition of $\varphi_1$ and
$\begin{pmatrix} R^t & 0\\0 &
I_{2q+2k}\end{pmatrix}:\Ze^{2p}\times 0^q\times
\Ze^{q+2k}\rightarrow \Ze^{2p}\times 0^q\times \Ze^{q+2k}$. Then
$\Ze^{2p}\times 0^q\times \Ze^{q+2k}=\triangle\oplus
\varphi(\Ze^n)$.
\end{lemma}
\begin{proof}
Let $y=\begin{pmatrix} y_1\\ 0 \\ y_2\\y_3 \end{pmatrix}$ in
$\Ze^{2p}\times 0^q\times \Ze^q\times \Ze^{2k}$ satisfying $y_3\in
n_1\Ze\times \cdots \times n_k\Ze\times n_1\Ze\times \cdots \times
n_k\Ze$.
Say $y_3=(n_1z_1 \cdots, n_kz_k, n_1z_{k+1}, \cdots,
n_kz_{2k})^t$. Then it is easy to see that $y\in \triangle$
exactly if
\begin{eqnarray*}
(R^t)^{-1}y_1&=&( m_1z_1, \cdots, m_kz_k,z_{k+1}, \cdots, z_{2k},
0,\cdots, 0)^t \quad \mbox{ and } \quad  y_2=0.
\end{eqnarray*}
It is clear from this that $\Ze^{2p}\times 0^q\times
\Ze^{q+2k}=\triangle\oplus \varphi(\Ze^n)$.
\end{proof}
Back to the proof of Proposition~\ref{embedding:prop}.
Putting $\varphi: \Ze^n\to (\Ze^{2p}\times 0^q\times \Ze^{q+2k})$
and $(\bar{T}^tJ)^{-1}:\Ze^{n+q+2k}\to H^*$ together, we get a map
$S:=(\bar{T}^tJ)^{-1}\circ \varphi:\Ze^n\to H^*$ with
$S(\Ze^n)=\mathcal{D}^{\perp}$. Let
\begin{eqnarray*}
Q_1=diag(d_1, \cdots, d_k),  & Q_2=diag(c_1, \cdots, c_k).
\end{eqnarray*}
A routine calculation shows that
\begin{eqnarray*}
S = \begin{pmatrix} W_1 & W_2 \\
\begin{pmatrix} 0 & -I_k & 0\\Q_2 & 0 & 0\end{pmatrix} & 0
\end{pmatrix},
\end{eqnarray*}
where $W_1$ and $W_2$ have sizes $(n+q)\times 2p$ and $(n+q)\times
q$ respectively:
\begin{eqnarray*}
W_1 =\begin{pmatrix}J_0(T^t_{11})^{-1}R^t\begin{pmatrix} P_1 & 0 & 0\\0 & P_1 & 0\\0 & 0 &  -I_{2p-2k}\end{pmatrix}\\
0\end{pmatrix}, &
W_2 =\begin{pmatrix} -J_0(T^t_{11})^{-1}T^t_{31} \\
I_q\\ T^t_{32}\end{pmatrix}.
\end{eqnarray*}

Clearly $S$ satisfies Definition~\ref{embedding map:def}(1)(2').
Then $S$ is an embedding map for
\begin{eqnarray*}
-\theta'=S^tJS = -\begin{pmatrix} \theta'_{11} & \theta'_{12}\\
\theta'_{21} & \theta'_{22}\end{pmatrix},
\end{eqnarray*}
where
\begin{eqnarray*}
\theta'_{11}&=& \begin{pmatrix} P_1 & 0 & 0\\0 & P_1 & 0\\0 & 0 &
-I_{2p-2k}\end{pmatrix}RF_{11}R^t\begin{pmatrix} P_1 & 0 & 0\\0 &
P_1 & 0\\0 & 0 &  -I_{2p-2k}\end{pmatrix} \\
&+&
\begin{pmatrix} 0 & -Q_2P_1 & 0\\ Q_2P_1 & 0 & 0 \\ 0 & 0 &
0\end{pmatrix},\\
\theta'_{12}&=&
\begin{pmatrix} P_1 & 0 & 0\\0 & P_1 & 0\\0 & 0 &
-I_{2p-2k}\end{pmatrix}RF_{11}\theta_{12},\\
\theta'_{21}&=& -\theta_{21}F_{11}R^t\begin{pmatrix} P_1 & 0 &
0\\0 & P_1 &
0\\0 & 0 & -I_{2p-2k}\end{pmatrix},\\
\theta'_{22}&=& -\theta_{21}F_{11}\theta_{12}+\theta_{22}.
\end{eqnarray*}

Proposition~\ref{module:prop} tells us that $\mathcal{S}(M)$ is a complete
Morita equivalence
$A^{\infty}_{\theta'}$-$A^{\infty}_{\theta}$-bimodule. Clearly the
dual $\phi^*:L^*_{\theta'}\rightarrow L^*_{\theta}$ of
$\phi:L_{\theta}\rightarrow L_{\theta'}$ is just
$-\tilde{T}^{-1}\circ \tilde{S}$. A routine calculation shows that
$\phi^*$ is given by the matrix
\begin{eqnarray*}
\mathscr{A} =-\begin{pmatrix}  -F_{11}R^t\begin{pmatrix} P_1 & 0 &
0\\0 & P_1 & 0\\0 & 0 &  -I_{2p-2k}\end{pmatrix}
 & -F_{11}\theta_{12}\\ 0 &
I_q\end{pmatrix}.
\end{eqnarray*}
It is also easy to see that the matrix form of the normalized
curvature $\frac{1}{2\pi i}\Omega$ is
\begin{eqnarray*}
\Phi=\begin{pmatrix} F_{11} & 0\\ 0 & 0 \end{pmatrix}.
\end{eqnarray*}
Now that we have the matrices $\theta, \theta', \mathscr{A}$, and $\Omega$,
Schwarz \cite[page 733]{Schwarz98} has shown how to find $g'=\begin{pmatrix} A' & B' \\
C' &  D' \end{pmatrix} \in SO(n,n|\Ze)$ such that
$\theta'=g'\theta$.  Actually we have the formulas:
\begin{eqnarray} \label{formula:eq}
C' = \mathscr{A}^{-1}\Phi, & D'= \mathscr{A}^{-1}-C'\theta, \\
\nonumber A'= \mathscr{A}^t+\theta'C', & B'=
\theta'\mathscr{A}^{-1}-A'\theta.
\end{eqnarray}
Our formulas (\ref{formula:eq}) are exactly the equation (53) of
\cite{Schwarz98}, in slightly different form. Straight-forward
calculations yield
\begin{eqnarray*}
C' =\begin{pmatrix}  \begin{pmatrix} T_4 & 0\\0 &
-I_{2p-2k}\end{pmatrix}(R^t)^{-1}
 & 0\\ 0 & 0\end{pmatrix}, &
D'= \begin{pmatrix}  \begin{pmatrix} 0 & -P_2 & 0\\P_2 & 0 & 0\\
0 & 0 & 0\end{pmatrix}R
 & 0\\ 0 & I_q\end{pmatrix}, \\
A'= \begin{pmatrix}  \begin{pmatrix} 0 & -Q_2 & 0\\Q_2 & 0 & 0\\0
& 0 &  0\end{pmatrix}(R^t)^{-1}
 & 0 \\ 0 & I_q\end{pmatrix}, &
B'= \begin{pmatrix}  \begin{pmatrix} Q_1  & 0 & 0\\0 & Q_1 & 0\\
0 & 0 & -I_{2p-2k}\end{pmatrix}R
 & 0 \\ 0 & 0 \end{pmatrix}.
\end{eqnarray*}
Let
\begin{eqnarray*}
\tilde{g}=
\begin{pmatrix} \tilde{A} & \tilde{B}\\ \tilde{C} & \tilde{D} \end{pmatrix}:=
g(g')^{-1}.
\end{eqnarray*}
Then $\tilde{g}\in SO(n,n|\Ze)$. A routine calculation shows that
$\tilde{C}=0$. By (\ref{O(n,n|R):eq})
we have $I=\tilde{A}^t\tilde{D}+\tilde{C}^t\tilde{B}=\tilde{A}^t\tilde{D}$.
Then $\tilde{A}$ is invertible.
Recall the matrix $\rho(\tilde{A})$ in Notation~\ref{rho:notation}.
We get
\begin{eqnarray*}
\tilde{g}=\begin{pmatrix} \tilde{A} & \tilde{B}\\ 0 & \tilde{D}
\end{pmatrix}
=\begin{pmatrix} I & \tilde{B}\tilde{A}^t\\ 0 & I
\end{pmatrix}\rho(\tilde{A}).
\end{eqnarray*}
Hence $\begin{pmatrix} I & \tilde{B}\tilde{A}^t\\ 0 & I
\end{pmatrix}=\tilde{g}(\rho(\tilde{A}))^{-1}\in SO(n, n|\Ze)$.
By (\ref{O(n,n|R):eq}) the matrix
$I^t(\tilde{B}\tilde{A}^t)=\tilde{B}\tilde{A}^t$ is
skew-symmetric. So we get
\begin{eqnarray*}
\tilde{g}=\nu(\tilde{B}\tilde{A}^t)\rho(\tilde{A}), & g=
\tilde{g}g'=\nu(\tilde{B}\tilde{A}^t)\rho(\tilde{A})g'.
\end{eqnarray*}
Notice that $g',\tilde{B}\tilde{A}^t$ and $\tilde{A}$ do not depend
on $\theta$. This finishes the proof of
Proposition~\ref{embedding:prop}.
\end{proof}

\begin{remark} \label{embedding:remark}
(1) We would like to point out that the argument right after (52)
of \cite{Schwarz98} is not quite complete. When $n=2$ the fact
that $W$ transforms the integral lattice $\wedge^{ev}(D^*)$ into
itself does not imply that $\tilde{W}$ has integral entries. In
other words, (\ref{formula:eq}) above may not give integral
matrices when $n=2$. So the case $n=2$ in \cite{Schwarz98} has to
be dealt separately, and it does follow from \cite{Rieffel81}. In
our situation we don't need to separate the case $n=2$ since $g'$
obviously has integral entries.

(2) Given $g'$ explicitly, one can also check directly that $g'
\in SO(n, n|\Ze)$ and $\theta'=g'\theta$: straight-forward
calculations show that $g'$ satisfies (\ref{O(n,n|R):eq}) and
$\theta'=g'\theta$. Then $g',\, \tilde{g}\in O(n, n|\Ze)$ and
hence we still have $I=\tilde{A}^t\tilde{D}$. Thus
$det(g')=det(\tilde{g}^{-1})=1$. Therefore $g'\in SO(n, n|\Ze)$.
\end{remark}

\begin{proof}[Proof of Theorem~\ref{orbit implies Morita
equiv:theorem}] We may think of $A_{\theta}$ as the universal
$C^*$-algebra generated by unitaries $\{U_{x, \theta}\}_{x\in
\Ze}$ satisfying the relation $U_{x, \theta}U_{y, \theta}=e(x\cdot
\theta y)U_{y, \theta}U_{x, \theta}$. For any $R\in GL(n|\Ze)$ and
$\theta\in \mathcal{T}_n$ there is a natural isomorphism
$A^{\infty}_{\theta}\rightarrow
A^{\infty}_{\rho(R)\theta}=A^{\infty}_{R\theta R^t}$ given by
$U_{x, \theta}\mapsto U_{(R^{-1})^tx,\rho(R)\theta}$. Under this
isomorphism $\delta_{X, \theta}$ becomes $\delta_{RX,
\rho(R)\theta}$ for any $X\in L^*$. Similarly, for any $N\in
\mathcal{T}_n\cap M_n(\Ze)$ and $\theta\in \mathcal{T}_n$ there is
a natural isomorphism $A^{\infty}_{\theta}\rightarrow
A^{\infty}_{\mu(N)\theta}=A^{\infty}_{\theta +N}$ given by $U_{x,
\theta}\mapsto U_{x,\mu(N)\theta}$, under which $\delta_{X,
\theta}$ becomes $\delta_{X, \mu(N)\theta}$ for any $X\in L^*$.
Now Theorem~\ref{orbit implies Morita equiv:theorem} follows from
Lemma~\ref{normalize:lemma} and Propositions~\ref{embedding:prop}
and \ref{module:prop}.
\end{proof}

\end{document}